\documentclass{acm_proc_article-sp}

\newdef{definition}{Definition}
\newtheorem{theorem}{Theorem}
\newtheorem{lemma}{Lemma}

\def\P{{\mathbb P}}
\def\E{{\mathbb E}}
\def\T{{\mathbb T}}
\def\D{{\cal D}}

\begin{document}
\title{Determining factors behind the PageRank log-log plot}

\numberofauthors{3}
\author{
\alignauthor
Yana Volkovich\\
      \affaddr{University of Twente}\\
      \affaddr{Dept. of Applied Mathematics, P.O. Box 217, 7500 AE }\\
      \affaddr{Enschede, The Netherlands}\\
      \email{y.volkovich@ewi.utwente.nl}
\alignauthor
Nelly Litvak\titlenote{The work is supported by NWO Meervoud grant no.~632.002.401 }\\
      \affaddr{University of Twente}\\
      \affaddr{Dept. of Applied Mathematics, P.O. Box 217, 7500 AE }\\
      \affaddr{Enschede, The Netherlands}\\
      \email{n.litvak@ewi.utwente.nl}
\alignauthor
Debora Donato \\
 \affaddr{Yahoo! Research}\\
 \affaddr{Barcelona Ocata 1, 1st floor 08003}\\
 \affaddr{Barcelona Catalunya, Spain}\\
 \email{debora@yahoo-inc.com}
} \maketitle
\begin{abstract}
We study the relation between PageRank and other parameters of
information networks such as in-degree, out-degree, and the
fraction of dangling nodes. We model this relation through a
stochastic equation inspired by the original definition of
PageRank. Further, we use the theory of regular variation to prove
that PageRank and in-degree follow power laws with the same
exponent. The difference between these two power laws is in a
multiple coefficient, which depends mainly on the fraction of
dangling nodes, average in-degree, the power law exponent, and
damping factor. The out-degree distribution has a minor effect,
which we explicitly quantify. Our theoretical predictions show a
good agreement with experimental data on three different samples
of the Web.
\end{abstract}




\subsection*{MSC 2000} 90B15, 68P10, 60J80


\section{Introduction}

Originally created for Web ranking, {\em PageRank} has become a
major method for evaluating popularity of nodes in information
networks. Besides its primary application in search engines,
PageRank is successfully used for solving other important problems
such as spam detection~\cite{Gyongyi04}, graph
partitioning~\cite{Andersen06}, and finding gems in scientific
citations~\cite{Chen06}, just to name a few.  The
PageRank~\cite{Brin98} is defined as a stationary distribution of a
random walk on a set of Web pages. At each step, with probability
$c$, the random walk follows a randomly chosen outgoing link, and
with probability $1-c$, the walk starts afresh from a page chosen at
random according to some distribution $f$. Such random jump also
occurs if a page is {\em dangling}, i.e. it does not have outgoing
links. In the original definition, the teleportation distribution
$f$ is uniform over all Web pages. Then the PageRank values satisfy
the equation
\begin{equation}
PR(i) = c \sum_{j \to i} \frac{1}{d_{j}} PR(j) + \frac{c}{n}
\sum_{j\in \D}PR(j)+\frac{1-c}{n}, \; i=1,\ldots,n,
\label{pr_brin_page_rw}
\end{equation}
where $PR(i)$ is the PageRank of page $i$, $d_j$ is the number of
outgoing links of page $j$, the sum is taken over all pages $j$ that
link to page $i$, $\D$ is a set of dangling nodes, $n$ is the number
of pages in the Web, and $c$ is the damping factor, which is a
constant between 0 and 1.

From equation (\ref{pr_brin_page_rw}) it is clear that the PageRank
of a page depends on popularity and the number of pages that link to
it. Thus, it can be expected that the distribution
 of PageRank should be related to the distribution of \emph{in-degree}, the number of
incoming links. Most of experimental studies of the Web agree that
in-degree follows a power law with exponent $\alpha=1.1$ for
cumulative plot, which corresponds to the famous value $2.1$ for the
density. Pandurangan et al.~\cite{Pandurangan02} discovered that
PageRank also follows a power law with the same exponent. Further
experiments~\cite{Becchetti06,Donato04,Fortunato05a} confirmed this
phenomenon. Mathematical justifications have been proposed
in~\cite{Avrachenkov06,Fortunato06} for the preferential attachment
models~\cite{Albert99}, and in~\cite{Litvak06}, where the relation
between PageRank and in-degree is modeled through a stochastic
equation.

At this point, it is important to realize that PageRank is a {\it
global} characteristic of the Web, which depends on in-degrees,
{\em out-degrees}, correlations, and other characteristics of the
underlying graph.  In contrast to in-degrees, whose impact on the
PageRank log-log plot is thoroughly explored and relatively well
understood, the influence of out-degrees and dangling nodes has
hardly received any attention in the literature. It is however a
common belief that dangling nodes are important~\cite{Eiron04}
whereas out-degrees (almost) do not affect the
PageRank~\cite{Fortunato05a}. We also note that in the literature,
there is no common agreement on the out-degree distribution. On
the Web data, Broder et al.~\cite{Broder00} report a power law
with exponent about 2.6 for the density, whereas e.g. Donato et
al.~\cite{Donato04} obtain a distribution, which is clearly not a
power law. On the other hand, for Wikipedia~\cite{Capocci06},
out-degree seems to follow a power law with the same exponent as
in-degree.

In the present paper we investigate the relations between PageRank
and in/out-degrees, both analytically and experimentally. Our
analytical model is an extension of~\cite{Litvak06}. We view the
PageRank of a random page as a random variable $R$ that depends on
other factors through a stochastic equation
resembling~(\ref{pr_brin_page_rw}).

It is clear that the PageRank values in~(\ref{pr_brin_page_rw})
scale as $1/n$ with the number of pages. In the analysis, it is more
convenient to deal with corresponding {\it scale-free} PageRank
scores
\begin{equation}
\label{eq_ri} R(i)=n PR(i),\quad i=1,\ldots,n,
\end{equation}
 assuming that $n$ goes to infinity.
In this setting, it is easier to compare the probabilistic
properties of PageRank and in/out-degrees, which are also
scale-free. In the remainder of the paper, by PageRank we mean the
scale-free PageRank scores~(\ref{eq_ri}). Then the original
definition~(\ref{pr_brin_page_rw}) can be written as
\begin{equation}
R(i) = c \sum_{j \to i} \frac{1}{d_{j}} R(j) + \frac{c}{n}\sum_{j\in
\D}R(j)+1-c, \; i=1,\ldots,n. \label{pr_brin_page}
\end{equation}
We are concerned with the {\it tail} probability $\P(R> x)$, i.e.
the fraction of pages with PageRank greater than $x$, when $x$ is
large. Our goal is to determine the asymptotic behavior of
$\P(R>x)$, that is, we want to find a known function $r(x)$ such
that $\P(R>x)/r(x)\to 1$ as $x\to\infty$. In this case, we say that
$\P(R>x)$ and $r(x)$ are asymptotically equivalent, which
essentially means that for large enough $x$, $\P(R>x)$ and $r(x)$
are close, and their log-log plots look the same.  We formally
describe power laws in terms of regular varying random variables,
and we use recent results on regular variation to obtain the
PageRank asymptotics. To this end, we provide a recurrent stochastic
model for the power iteration algorithm commonly used in PageRank
computations~\cite{Langville03}, and we obtain the PageRank
asymptotics after each iteration.

The analytical results suggest that the PageRank and in-degree
follow power laws with the same exponent. The out-degrees and
dangling nodes affect only a multiple factor, for which we find an
exact expression. It follows that the out-degree sequence has a
truly minor influence whereas the fraction of dangling nodes has a
slightly greater impact on the multiple coefficient. The
experiments on the Indochina-2004 Web sample~\cite{IndoChina}, on
the EU-2005 Web sample~\cite{IndoChina}, and on the Stanford
Web~\cite{Stanford}, show that our model correctly predicts the
evolution of the PageRank distribution through the series of power
iterations, and it adequately captures the influence of the
network parameters.

\section{Preliminaries}\label{section_regvar}

We start with preliminaries on the theory of regular variation,
which is a natural formalization of power laws. More comprehensive
details could be found, for instance, in~\cite{Bingham_RV}. We
also refer to Jessen and Mikosch~\cite{Jessen06} for an excellent
recent review.

\begin{definition} A function
$L(x)$ is {\em slowly varying} if for every $t>0$,
\begin{equation*}
\frac{L(tx)}{L(x)}\to 1\quad\mbox{ as }\quad x\to \infty.
\end{equation*}
\end{definition}
\begin{definition}
A non-negative random variable $X$ is said to be {\em regularly
varying} with index $\alpha$ if
\begin{equation}
\P(X>x)\sim x^{-\alpha}L(x) \quad\mbox{ as }\quad x\to
\infty,\label{e02}
\end{equation}
for some positive slowly varying function $L(x)$.
\end{definition}

Here, as in the remainder of this paper, the notation $a(x)\sim
b(x)$ means that $a(x)/b(x) \to 1.$

The {\em asymptotic} equivalence (\ref{e02}) is a formalization of
a power law. In words, it means that for large enough $x$, the
{tail} distribution $\P(X>x)$ can be approximated by the
regularly varying function $x^{-\alpha}L(x)$, which is, in turn,
approximately proportional to $x^{-\alpha}$ due to the definition
of $L$.

Regularly varying random variables represent a subclass of a much broader class of long-tailed random variables.

\begin{definition}
A  random variable $X$ is  {\em
long-tailed} if for any $y>0$,
\begin{equation}
\P(X>x+y)\sim \P(X>x)\quad\mbox{ as }\quad x\to
\infty\label{eq_long}.
\end{equation}
\end{definition}

Next lemma describes the behavior of a product and random sums of
regular varying random variables. The relation (i) is known as
Breiman's theorem (see e.g. Lemma 4.2.(1) in~\cite{Jessen06}).
Properties (ii) and (iii) are, respectively, statements (2) and
(5) of Lemma 3.7 in~\cite{Jessen06}.

\begin{lemma}\label{lemma properties}
\begin{description}
\item[(i)] Assume that $X_1$ and $X_2$ are two independent
non-negative random variables such that $X_1$ is regularly varying
with index $\alpha$ and that
$\mathbb{E}(X_2^{\alpha+\epsilon})<\infty$ for some $\epsilon>0.$
Then
\begin{equation*}
    \P(X_1 X_2>x)\sim\mathbb{E}(X_2^{\alpha})\P(X_1>x).
\end{equation*}

\item[(ii)]Assume that $N$ is regularly varying with index
$\alpha\geq 0;$ if $\alpha=1$, then assume that
$\mathbb{E}(N)<\infty.$ Moreover, let $(X_{i})$ be i.i.d. sequence
such that $\mathbb{E}(X_1)<\infty$ and $\P(X_1>x)=o(\P(N>x)).$ Then
as $x\to\infty,$
\begin{equation*}
    \P\left(\sum_{i=1}^{N}X_{i}>x\right)\sim(\mathbb{E}(X_1))^{\alpha}\P(N>x).
\end{equation*}

\item[(iii)] Assume that $\P(N>x)\sim r \P(X_1>x)$ for some $r>0$, that
$X_1$ is regularly varying with index $\alpha\geq 1$, and
$\E(X_1)<\infty.$ Then
\begin{equation*}
\P\left(\sum_{i=1}^{N}X_{i}>x\right)\sim
(\mathbb{E}(N)+r(\mathbb{E}(X_1))^{\alpha})\P(X_1>x).
\end{equation*}
\end{description}
\end{lemma}

\section{The Model}\label{section_model}
\subsection{In-degree}\label{ssection_inout}

It is a common knowledge that in-degrees in the Web graph obey a
power law with exponent about $2.1$ for the density, which
corresponds to $1.1$ for cumulative plot. The power law exponent may
deviate somewhat depending on a data set~\cite{Baeza06} and an
estimator~\cite{Newman05}. As in our previous work~\cite{Litvak06},
we model the in-degree as an integer regularly varying random
variable. To this end, we assume that the in-degree of a random page
is distributed as $N(T)$, where $T$ is regularly varying with index
$\alpha$ and $N(t)$ is the number of Poisson arrivals on the time
interval $[0,t]$, when arrival rate is 1. If $T$ is regularly
varying then $N(T)$ is also regularly varying and asymptotically
identical to $T$ (see e.g.~\cite{Litvak06}). Thus, $N(T)$ is indeed
integer and obeys the power law. To simplify the notation, we will
use $N$ instead of $N(T)$ throughout the paper. The proposed
formalization for the in-degree distribution allows us to model the
number of terms in the summation in (\ref{pr_brin_page}).

\subsection{Out-degree and inspection paradox}

 Now, we want to model the weights $1/d_j$ in~(\ref{pr_brin_page}). Recall that $d_j$ is
the out-degree of page $j$ that has a link to page $i$. In
\cite{Litvak06} we studied the relation between in-degree and
PageRank assuming that out-degrees of all pages are constant, equal
to the expected in-degree $d$. In this work, we make a step further
allowing for random out-degrees.

We model out-degrees of pages linking to a randomly chosen page as
independent and identically distributed random variables with
arbitrary distribution. Thus, consider a random variable $D$,
which represents the out-degree of a page that links to a
particular randomly chosen page $i$. Note that $D$ is \emph{not}
the same random variable as an out-degree of a random page since
the additional information that a page has a link to $i$, alters
the out-degree distribution. This famous phenomenon, called
\emph{inspection paradox}, finds its mathematical explanations in
Renewal Theory. The inspection paradox roughly states that an
interval containing a random point tends to be larger than a
randomly chosen interval~\cite{Ross_SP}. For instance, in~\cite{Ross03}, a number of children in a family,
to which a randomly chosen child belongs, is stochastically larger
than a number of children in a randomly chosen family. Likewise, a
number of out-links $D$ from a page containing a random link,
should be stochastically larger than an out-degree of a random
page. We will refer to $D$ as \emph{effective out-degree}. The
term is motivated by the fact that the distribution of $D$ {is}
{the one} that participates in the PageRank formula.

Now, let $p_j$ be a fraction of pages with out-degree $j\ge 0$. Then
we have
\begin{equation}
\label{eq_D}
    \lim_{n\to\infty}\P(D=j)=\frac{jp_j}{d},\quad j\ge 1.
\end{equation}
where $d$ is the average in/out-degree, and $n$ is the number of
pages in the Web. For sufficiently large networks, we may assume
that the distribution of $D$ equals to its limiting distribution
defined by (\ref{eq_D}). Note that, naturally, the probability that
a random link comes from a page with out-degree $j$ is proportional
to $j$. This was implicitly observed by Fortunato et al.
in~\cite{Fortunato05a}, who in fact used (\ref{eq_D}) in their
computations for the mean-filed approximation of PageRank.

\subsection{Stochastic equation}\label{ssection_steq}

We view the scale-free PageRank of a random page as a random
variable $R$ with $\E(R)=1$. Further, we assume that the PageRank of
a random page does not depend on the fact whether the page is
dangling. Indeed, it can be shown that the PageRank of a page can
not be altered significantly by modifying outgoing
links~\cite{Avrachenkov06a}. Moreover, experiments e.g.
in~\cite{Eiron04} show that dangling nodes are often just regular
pages whose links have not been crawled, for instance, because it
was not allowed
 by \verb"robot.txt". Besides, even authentically dangling
pages such as \verb".pdf" or \verb".ps" files, often contain
important information and gain a high ranking independently of the
fact that they do not have outgoing links. We note that such
independence implies that the average PageRank of dangling nodes
is~1, and thus the fraction of the total PageRank mass concentrated
in dangling nodes, equals to the fraction of dangling nodes $p_0$:
\begin{equation*}
p_0=\frac{1}{n}\sum_{j\in \D}R(j).
\end{equation*}

Our goal is to model and analyze to which extent the tail
probability $\P(R>x)$ for large enough $x$ depends on the
in-degree $N$, the effective out-degree $D$, and the fraction of
dangling nodes $p_0$. To this end, we model PageRank $R$ as a
solution of a stochastic equation involving $N$ and $D$.  Inspired
by the original formula (\ref{pr_brin_page}), the stochastic equation for
the scale-free PageRank is as follows:
\begin{equation}\label{main_st_eq}
    R\stackrel{d}=c\sum_{j=1}^{N}\frac{1}{D_{j}}R_{j}+[1-c(1-p_0)].
\end{equation}

 Here $N$, $R_j$'s and $D_j$'s are independent; $R_j$'s
are distributed as $R$,  $D_j$'s are distributed as $D$, and
$a\stackrel{d}{=}b$ means that $a$ and $b$ have the same probability
distribution. As before, $c\in(0,1)$ is a damping factor.

We note that the independence assumption for PageRanks and
effective out-degrees of pages linking to the same page, is
obviously not true in general. However, there is also no direct
relation between these values as there is no experimental evidence
that such dependencies would crucially influence the PageRank
distribution. Thus, we assume independence in this study.

The stochastic equation (\ref{main_st_eq}) is a generalization of
the equation analyzed in~\cite{Litvak06},  where it was assumed that
$D_j$'s are constant. In order to demonstrate applicability of our
model, we will use (\ref{main_st_eq}) to derive a mean-field
approximation for the PageRank of a page with given in-degree. It
follows from~(\ref{eq_D}) that
\begin{equation*}
\label{eq_1d}
    \E\left(\frac{1}{D}\right)=\sum_{j=1}^{\infty}\frac{1}{j}\P(D=j)=
    \sum_{j=1}^{\infty}\frac{1}{j}\,\frac{jp_j}{d}=\frac{1-p_0}{d}.
\end{equation*}
Then, assuming that $\E(R_j)=1$, $j=1,2,\ldots$,  we obtain
\begin{equation}
\label{eq_mean} \E(R|N) =\frac{c(1-p_0)}{d}\,N+[1-c(1-p_0)].
\end{equation}
If $p_0=0$ then this coincides with the mean-field approximation
by Fortunato et al. in~\cite{Fortunato05a}, obtained directly from
the PageRank definition under minimal independence assumptions and
without considering dangling nodes.

Equation (\ref{main_st_eq}) belongs to the class of stochastic
recursive equations that were discussed in detail in the recent
survey by Aldous and Bandyopadhyay~\cite{Aldous05}. In particular,
(\ref{main_st_eq}) has an apparent similarity with distributional
equations motivated by branching processes and branching random
walks. Such equations were studied in detail by Liu in \cite{Liu01}
and his other papers. Taking expectations in (\ref{eq_mean}), we see
that if $\E(R_j)=1$, $j=1,2,\ldots$, then $\E(R)$ also equals 1. In
Section~\ref{section_analysis} we will show that (\ref{main_st_eq})
has a unique solution $R$ such that $\E(R)=1$.

\section{Model for power iterations}\label{section_pwit}

In this section, we will introduce an iteration procedure for
solving (\ref{main_st_eq}). This procedure can be seen as a
stochastic model for the power iteration method commonly used in
PageRank computations. We first present the notations, which are in
lines with Liu~\cite{Liu01}.

Let
$\left\{\left(N_{u},\frac{1}{D_{u_{1}}},\frac{1}{D_{u_{2}}},\ldots\right)\right\}_{u}$
be a family of independent copies of
$\left(N,\frac{1}{D_{1}},\frac{1}{D_{2}},\ldots\right)$ indexed by
all finite sequences\\ $u=u_{1}\ldots u_{n},\mbox{ }u_{i}\in
\{1,2,\ldots\}$. And let $\mathbb{T}$ be the Galton-Watson tree with
defining elements $\{N_{u}\}:$ we have $\emptyset\in\mathbb{T}$ and,
if $u\in\mathbb{T}$ and $i\in\{1,2,\ldots\}$, then concatenation
$ui\in \mathbb{T}$ if and only if $1\leq i \leq N_{u}.$ In other
words, we indexed the nodes of the tree with root $\emptyset$ and
the first level nodes $1,2,..N_{\emptyset},$ and at every subsequent
level, the $i$th offspring of $u$ is named $ui$ (see
Figure~1).{\sloppy

}
\begin{figure}\label{gwtree_example} \centering
\psfig{file=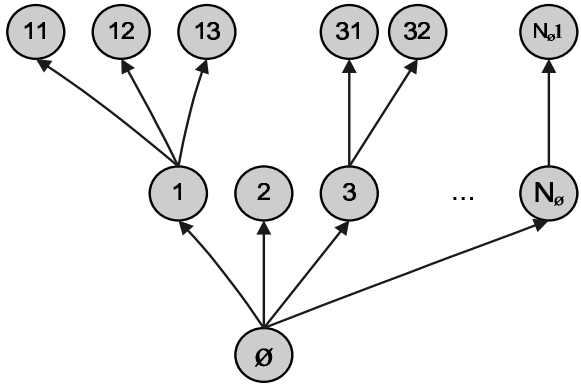, height=2in, width=3in,} \caption{An example
of Galton-Watson tree}
\end{figure}

Now, we will iterate the equation (\ref{main_st_eq}). We start with
initial distribution $R^{(0)}$, $\E\left(R^{(0)}\right)=1$, and for
every $k\geq 1$, we define the result of the $k$th iteration through
a distributional identity
\begin{equation}\label{eq_iter}
R^{(k)}\stackrel{d}{=}c\sum_{j=1}^{N}\frac{1}{D_{j}}R^{(k-1)}_{j}+[1-c(1-p_0)],
\end{equation}
where $N$, $R_j^{(k-1)}$ and $D_j$, $j\ge 1$, are independent. We
argue that if $R^{(0)}\equiv1$ then $R^{(k)}$ serves as a stochastic
model for the result of the $k$th power iteration in standard
PageRank computations. Indeed, according to (\ref{eq_iter}) for
$R^{(1)}$ we can obtain
\[R^{(1)}\stackrel{d}{=}c\sum_{j=1}^{N}\frac{1}{D_{j}}+[1-c(1-p_0)],\]
which clearly corresponds to the first power iteration with
initial uniform vector:
\[PR^{(1)}(i)=c\sum_{j\to i}\frac{1}{d_{j}}+[1-c(1-p_0)],\; i=1\ldots n.\]
This argument can be easily extended to further iterations.

Since PageRank vector is always a result of a finite number of
iterations, it follows that $R^{(k)}$ describes the distribution
of PageRank if the power iteration algorithm stops after $k$
steps. Assuming that in-degrees, effective out-degrees and
$R^{(0)}_u$, $u\in\T$, are independent, and repeatedly applying
(\ref{eq_iter}), we derive the following representation for
$R^{(k)}$:
\begin{align} \nonumber
   R^{(k)}=c^{k}\sum_{u=u_{1}..u_{k}\in\mathbb{T}}\frac{1}{D_{u_1}}\ldots\frac{1}{D_{u_1..u_k}}R^{(0)}_{u_1..u_k}&\\
\label{iter_long}  +
[1-c(1-p_0)]\sum_{n=0}^{k-1}c^{n}Y^{(n)},\quad k\ge 1,&
\end{align}
where
\begin{equation*}\label{eq_y}
Y^{(n)}=\sum_{u=u_{1}\ldots
u_{n}\in\mathbb{T}}\frac{1}{D_{u_{1}}}\ldots\frac{1}{D_{u_{1}\ldots
u_{n}}},\;n\ge 1.
\end{equation*}
The random variable $Y^{(n)}$ represents the sum of the weights of
the $n$th level of the Galton-Watson tree, where the root has
weight $1$, each edge has a random weight distributed as $1/D$,
and the weight of a node is a product of weights of the edges,
which are on the way from the root to this node.

In the subsequent analysis we will prove that iterations
$R^{(k)}$, $k\ge 1$, converge to a unique solution of
(\ref{main_st_eq}), and we will obtain the tail behavior of
$R^{(k)}$ for each $k\ge 1$. This will give us the asymptotic
behavior of the PageRank vector after an arbitrary number of power
iterations.

\section{Analytical results}
\label{section_analysis}

First, we establish that our main stochastic equation
(\ref{main_st_eq}) indeed defines a unique distribution $R$, that
can serve as a model for the PageRank of a random page. The result
is formally stated in the next theorem (the proof is given in
Section~\ref{section_proofs}).

\begin{theorem}\label{th_solut}
Equation (\ref{main_st_eq}) has a unique non-trivial solution with
mean $1$ given by
\begin{equation}\label{eq_solut}
R^{(\infty)}=\lim_{k\to\infty}R^{(k)}=[1-c(1-p_0)]\sum_{n=0}^{\infty}c^{n}Y^{(n)}.
\end{equation}
\end{theorem}

Now we are ready to describe the tail behavior of $R^{(k)}$,
$k\ge 1$, which models the PageRank after $k$ power iterations. The main result is presented in
Theorem~\ref{th_rv} below.

\begin{theorem}
\label{th_rv} If $\P\left(R^{(0)}>x\right)=o(\P(N>x))$, then for all
$k\ge 1$,
    \begin{equation*}\label{eq_rv_iter}
    \P(R^{(k)}>x)\sim C_{k}\P(N>x)\mbox{ as
    }x\to \infty,
    \end{equation*}
    where
    $C_{k}=\left(\frac{c(1-p_0)}{d}\right)^{\alpha}\sum_{j=0}^{k-1}c^{j\alpha}b^{j}$, and $b=d\E\left(1/D^{\alpha}\right)=\sum_{j=1}^{\infty}\frac{p_j}{j^{\alpha-1}}$.
\end{theorem}

The form of the coefficient $C_k$ arises from the proof, which
relies on the results from~\cite{Jessen06}. The proof is provided
in Section~\ref{section_proofs}. For large enough $k$, $C_k$ can
be approximated by
\[C=\lim_{k\to\infty}C_k=\frac{c^{\alpha}(1-p_0)^{\alpha}}{d^\alpha(1-c^\alpha b)}.\]
From the Jensen's inequality $\E(1/D^\alpha)\ge (\E(1/D))^{\alpha}$ and (\ref{eq_1d}), it follows that $b\ge (1-p_0)^{\alpha} d^{1-\alpha}$, and hence,
\begin{equation}
\label{eq_c} C\ge
\frac{c^{\alpha}(1-p_0)^{\alpha}}{d^\alpha(1-c^\alpha
(1-p_0)^{\alpha} d^{1-\alpha})}.\end{equation} The last expression
is the value of $C$ if out-degree of all non-dangling nodes is  a
constant. Note that if $\alpha\approx 1.1$, then the difference
between the left- and the right-hand sides of (\ref{eq_c}) is really
small for any reasonable out-degree distribution.

From Theorem~\ref{th_rv} we can make interesting conclusions about
the relation between PageRank and in/out-degrees. As it is commonly
known from experiments, the power law exponent of the PageRank is
the same as the power law exponent of in-degree. Clearly,  this
exponent is not affected by out-degrees. Thus, in-degree remains a
major factor shaping the PageRank distribution. The multiple factor
$C_k$, $k\ge 1$, depends mainly on the mean in-degree $d$, damping
factor $c$, and the fraction of non-dangling nodes $(1-p_0)$. The
values $p_j$, $j\ge 1$, that specify the out-degree distribution,
have some effect on the coefficient $b$ but this results in a truly
minor impact on the PageRank asymptotics. Hence, our results confirm
the common idea that the out-degree distribution has a very little
influence on the PageRank, but here we could also explicitly
quantify this minor effect. In the next section we will compare out
analytical findings with experimental results.

\section{Experiments}\label{section_numres}

We performed experiments on Indochina-2004 and EU-2005 Web samples
collected by The Laboratory for Web Algorithmics (LAW),
Dipartimento di Scienze dell'Informazione (DSI) of the Università
degli studi di Milano~\cite{IndoChina}. We also used a
Stanford-2002 Web sample~\cite{Stanford}.  In Figures~2--4 below
we present cumulative log-log plots for in-degree/PageRank. The
$y$-axis corresponds to the fraction of pages with
in-degree/PageRank greater than the value on the $x$-axis. For
in-degree, the power law exponent in evaluated using the maximum
likelihood estimator from~\cite{Newman05}, and the straight line
is fitted accordingly. For the PageRank, we plot the {\em
theoretically predicted} straight lines obtained from
Theorem~\ref{th_rv}.

The Indochina set contains 7414866 nodes and 194109311 links. The
results are presented in Figure~2 below. The in-degree plot
resembles a power law except for the excessively large fraction of
pages with in-degree about $10^4$. We suspect that this irregularity
might be related to the specific crawling
technique~\cite{Becchetti06a}. For more detail on this data set
see~\cite{Baeza06}. For Indochina, we obtain a power law exponent
$1.17$ for cumulative plot, which is quite different from the result
in~\cite{Baeza06}. This demonstrates the sensitivity of estimators
for the power law exponent. Indeed, the exponent
$0.6$~in~\cite{Baeza06} reflects the behavior in the first part of
the plot, whereas $1.17$ gives more weight on the tail of the
in-degree distribution.{\sloppy

}

\begin{figure}[hbt]\label{pic_ind} \caption{Indochina data set: cumulative log-log plots
for in-degree/PageRank. The straight lines for the PageRank plots
are predicted by the model.}
            \centering {\epsfxsize=1.5in \epsfysize=1.1in \epsfbox{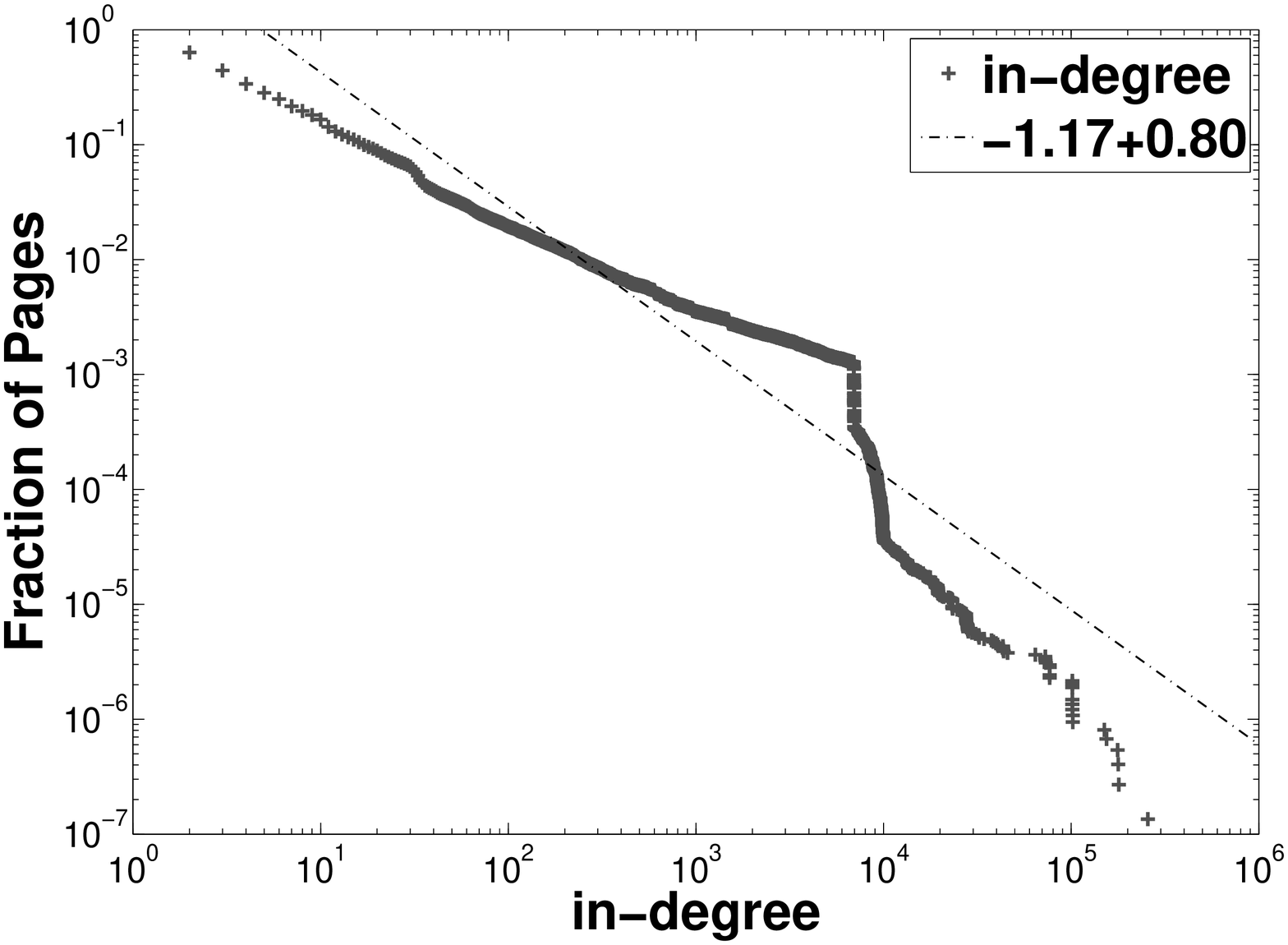}
            \epsfxsize=1.5in \epsfysize=1.1in \epsfbox{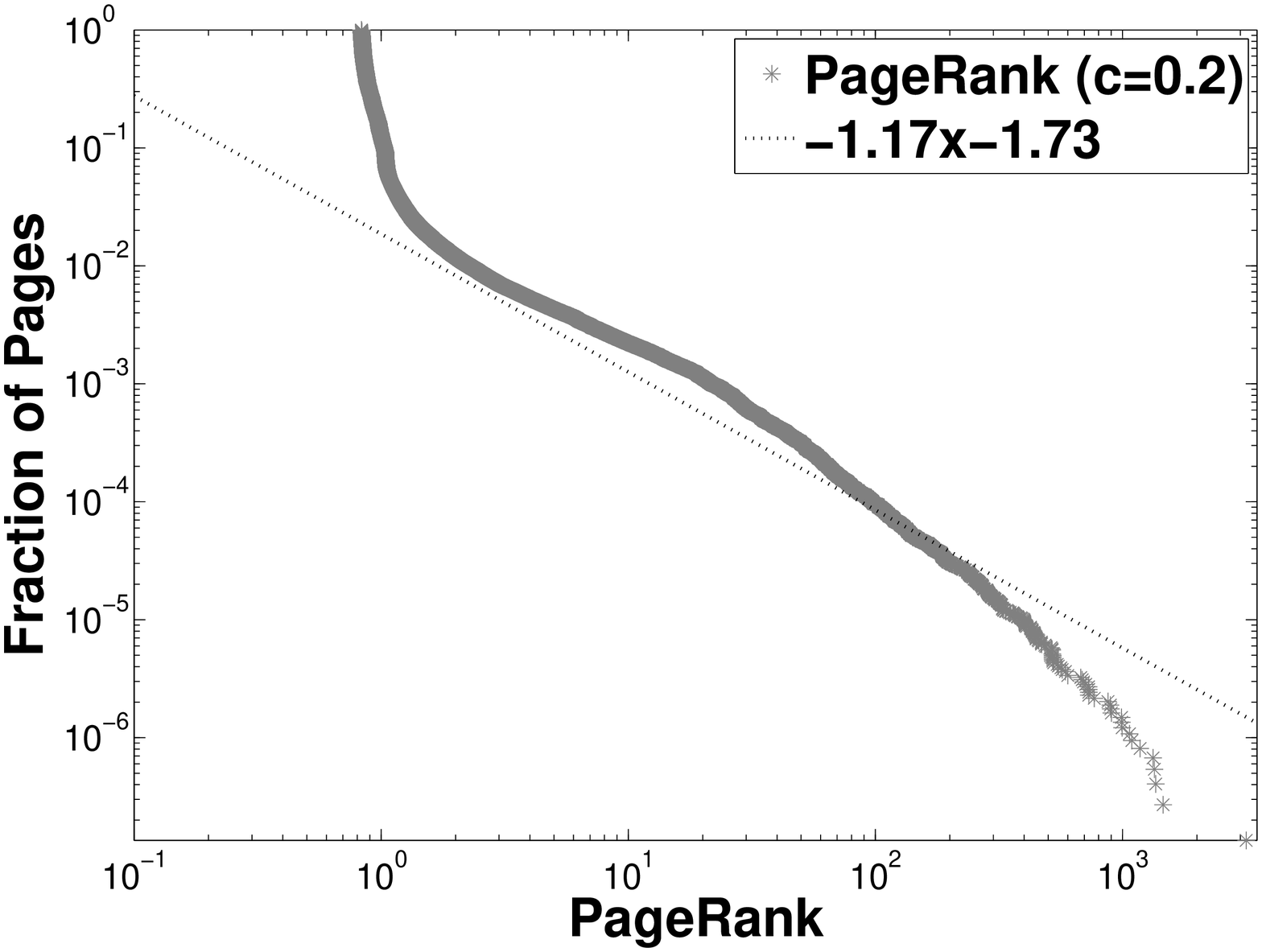}}
            \centering {\epsfxsize=1.5in \epsfysize=1.1in \epsfbox{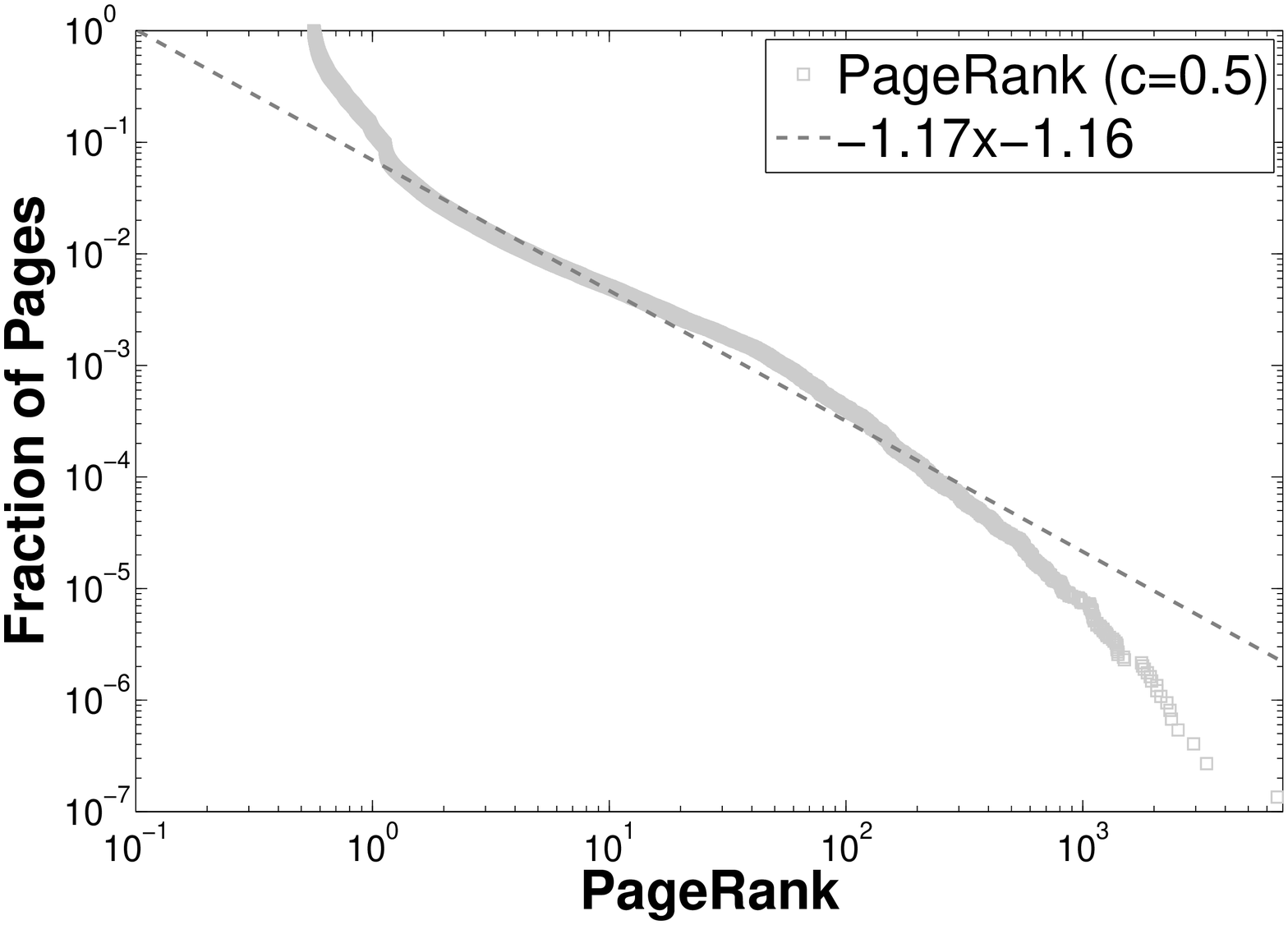}
            \epsfxsize=1.5in \epsfysize=1.1in \epsfbox{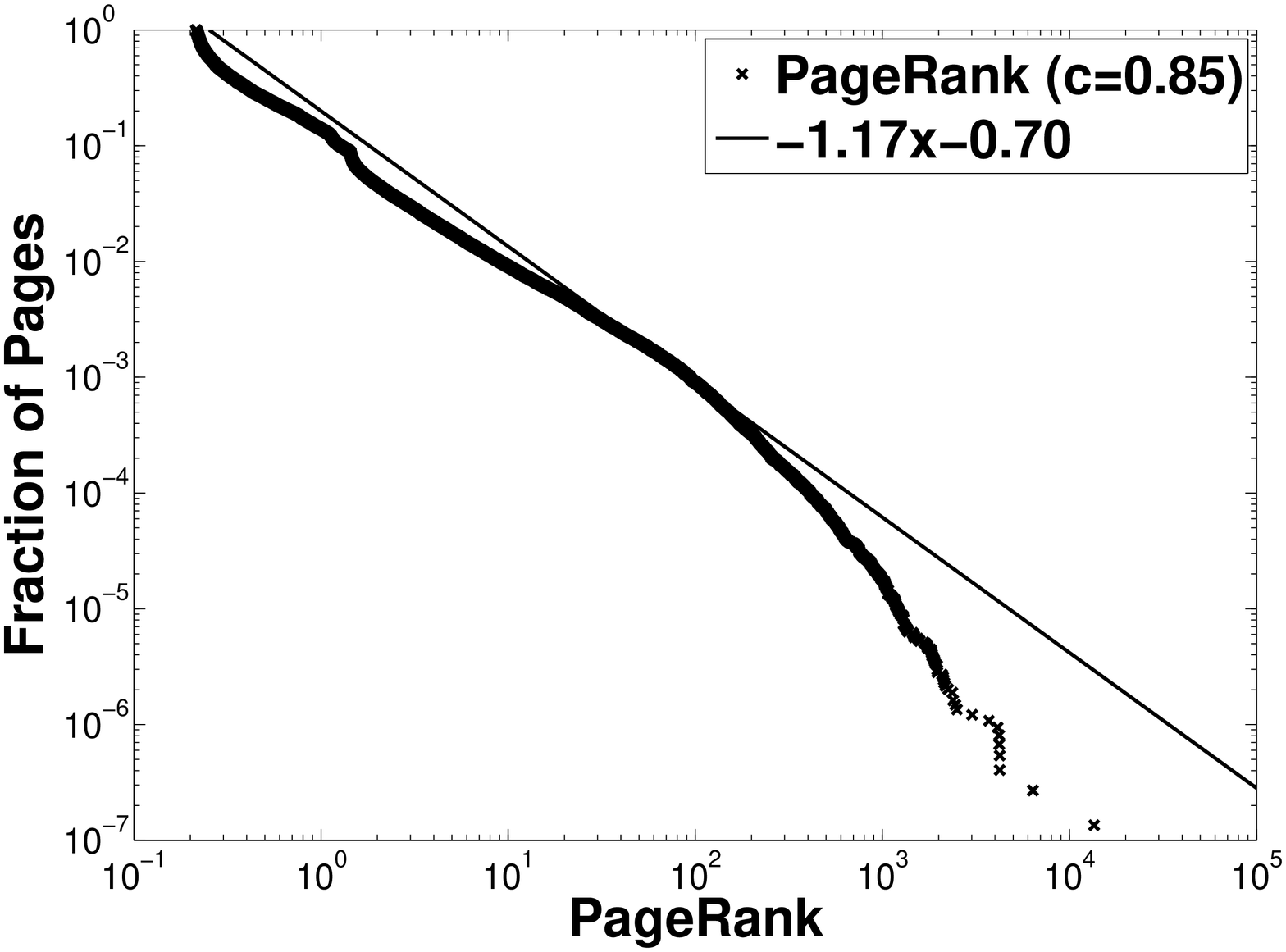}}
\end{figure}

We fit the straight line $y=-1.17x+0.80$ into the in-degree plot and
then compute the distance
\[\log_{10}(C)=\log_{10}\left(\frac{c^{\alpha}(1-p_0)^{\alpha}}{d^\alpha(1-c^\alpha b)}\right)\]
between the in-degree and the PageRank log-log plots for $c=0.2,
0.5$, and $0.85$. With $d=26.17$, $p_0=0.18$, and $b=0.65$, we
obtain the following prediction for the PageRank log-log plot:
$y=-1.17x-1.73$ for $c=0.2$, $y=-1.17x-1.16$ for $c=0.5$, and
$y=-1.17x-0.70$ for $c=0.85$. In Figure~\ref{pic_ind} we show these
{\it theoretically predicted} lines and the experimental PageRank
log-log plots. We see that for this data set, our model provides the
linear fit with a striking accuracy.{\sloppy

}

We performed the same experiment for EU-2005 of 862664 nodes and
19235140  links. In this data set in-degree shows a typical power
law behavior, which is fitted perfectly by $y=-1.1x+0.61$. We use
the same approach to calculate the difference between the in-degree
and PageRank plots for $d=22.3$, $p_0=0.08$, $b=0.70.$ Thus, the
theoretical prediction for the PageRank are $y=-1.1x-1.63$,
$y=-1.1x-1.07$, and $y=-1.1x-0.60$ for $c=0.2, 0.5$, and $0.85$,
respectively. The log-log plots for experimental data, the fitted
straight line for in-degree, and corresponding theoretical straight
lines for PageRank, are presented in Figure~3.
\begin{figure}[hbt]\label{pic_eu} \caption{EU-2005 data set: cumulative log-log plots
for in-degree/PageRank. The straight lines for the PageRank plots
are predicted by the model.}\centering {\epsfxsize=3.5in
\epsfysize=2.5in \epsfbox{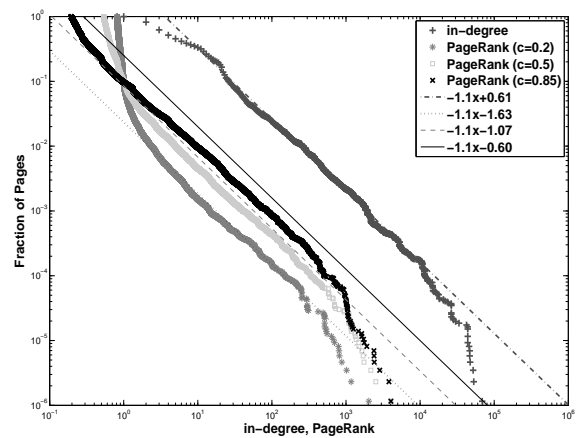}}

\end{figure}

Finally, we verify out model for power iterations. For that, we
use a smaller Web sample from~\cite{Stanford} that contains
$281903$ pages and above $2.3$ million links. In Figure~4
 we show the cumulative log-log plot of in-degree, and the log-log plots of the PageRank after the
$1$st, the $2$nd, and the last power iterations for the damping
factor $0.85$. To predict the difference between in-degree and
PageRank's iterations we use the result of Theorem~\ref{th_rv} for
$d=8.2032$, $p_0=0.006$, and $b=0.8558.$ Thus, if in-degree
distribution could be fitted by $y=-1.1x+0.08,$ then $y=-1.1x-1.00,$
$y=-1.1x-0.77,$ and $y=-1.1x-0.46.$ are the predicted PageRank after
the $1$st, the $2$nd, and the last power iterations, respectively.
Although the obtained lines do not match perfectly the PageRank
distribution, we see that our model correctly captures the dynamics
of the PageRank distribution in successive power iterations. The
difference between the theoretical prediction and the real data
might occur because of the specific structure of this data set. For
instance, the number of dangling nodes in this Web sample is
negligibly small, which is not true for the real Web.
\begin{figure}[hbt]\label{pic_stan} \caption{Stanford data set:
cumulative log-log plots for in-degree/PageRank. The straight
lines for the PageRank plots are predicted by the model for the
$1$st, the $2$nd, and the last power iterations.}
            \centering {\epsfxsize=3.5in \epsfysize=2.5in  \epsfbox{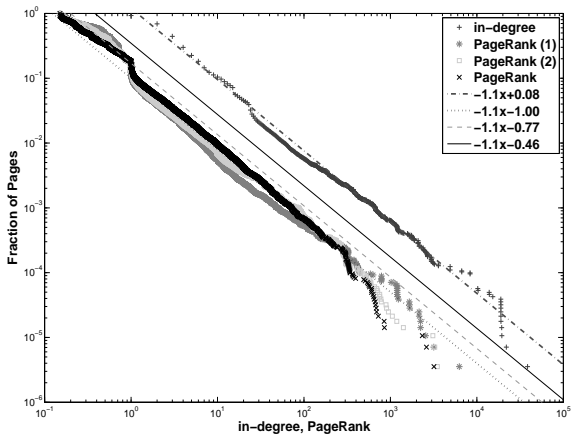}}
\end{figure}

\section{Discussion}\label{section_conclusions}

In this paper, we proposed an analytical stochastic model that
helps to predict the shape of the PageRank log-log plot on basis
of in-degree distribution, the damping factor, and the fraction of
dangling nodes. It also follows form the model that the out-degree
distribution has a truly minor impact on the PageRank. To make our
mathematical model analytically tractable, we had to allow for
several simplifying assumptions, such as independence of certain
parameters and uniform teleportation. Experiments show that our
theoretical model matches the Web data with a good accuracy.

One can argue that a uniform teleportation vector $f$ does not suit
anymore for Web ranking~\cite{Eiron04}. Indeed, there are smarter
choices of $f$ that take into account user's preferences, favor
certain topics related to a query~\cite{Haveliwala03}, or give
higher weights to trusted pages for eliminating the
spam~\cite{Eiron04}. The goal of this paper however was not
improving the Web ranking but rather analyzing why the PageRank
vector has certain properties reflected in its log-log plot. In
order to capture the influence of in- and out-degrees, we had to
make simplifying assumptions on other factors. However, we believe
that our approach is promising in modeling relations between
different parameters in the Web. In further research, we plan to
gradually improve our model including dependencies, personalization,
and other important factors relevant for the contemporary Web
search.

\section{Proofs}\label{section_proofs}\smallskip
\begin{proof}[of Theorem~\ref{th_solut}]
First, we establish that $R^{(\infty)}$ is well-defined random
variable. We consider some initial distribution $R^{(0)}$ with
$\E(R^{(0)})=1$. Then the first part of (\ref{iter_long}) has a mean
$c^k(1-p_0)^k$, and hence it converges in probability to 0 because,
by the Markov inequality, the probability that this term is greater
than some $\epsilon>0$ is at most $c^k(1-p_0)^k/\epsilon\to 0$ as
$k\to\infty$. Further, since $(1-p_0)^{-n}Y^{(n)}$ is a martingale
with mean~1, and $\lim_{n\to\infty}(1-p_0)^{-n}Y^{(n)}$ exists and
it is finite (see \cite{Liu01}), the second part of
(\ref{iter_long}) converges a.s. to $R^{(\infty)}$ as $k\to \infty$.
It follows that (\ref{iter_long}) converges to $R^{(\infty)}$ in
probability and according to the monotone convergence theorem
\begin{equation*}
\E\left(R^{(\infty)}\right)=[1-c(1-p_0)]\lim_{k\to\infty}\sum_{n=1}^{k}c^n\mathbb{E}\left(Y^{(n)}\right)=1.
\end{equation*}

It is easy to verify that $R^{(\infty)}$ in (\ref{eq_solut}) is a
solution of (\ref{main_st_eq}). To prove the uniqueness, we assume
that there is another solution with mean~1, then we take this
solution as an initial distribution $R^{(0)}$ and repeat the
argumentation above. Thus, we can conclude that there is no other
fixed point of (\ref{main_st_eq}) with mean~1 except $R^{(\infty)}$.
\end{proof}

\begin{proof}[of Theorem~\ref{th_rv}]
We will use the induction. For $k=1$, we derive
 \begin{align*}
   \P&\left(R^{(1)}>x\right)\sim\P\left(\sum_{j=1}^{N}\frac{c}{D_{j}}R_j^{(0)}+[1-c(1-p_0)]>x\right)\\
   &\sim
   \left(\frac{c(1-p_0)}{d}\right)^{\alpha}\P(N>x-[1-c(1-p_0)])\\
   &\sim C_1\P(N>x)\;\mbox{as}\;x\to\infty,
 \end{align*}
where the second relation follows from Lemma~\ref{lemma
properties}$(ii)$ because $\mathbb{E}(N)=d<\infty,$
$\E\left(R_1^{(0)}\right)=1,$
$\mathbb{E}\left({c}D_1^{-1}R_1^{(0)}\right)={c(1-p_0)}{d^{-1}}~<\infty$,
and $\P\left({c}D_1^{-1}R_1^{(0)}>x\right)=o(\P(N>x))$, and the last
relation follows from (\ref{eq_long}).

Now, assume that  the result has been shown for $(k-1)$th
iteration, $k\ge 2$. Then Lemma \ref{lemma properties}$(i)$ yields
\begin{align*}
    \P\left(\frac{c}{D}R^{(k-1)}>x\right)
    &\sim c^{\alpha}\mathbb{E}\left(\frac{1}{D^{\alpha}}\right)C_{k-1}\P(N>x)\\
    &=\frac{c^{\alpha}}{d}\;b\; C_{k-1}\P(N>x),
\end{align*}
where
\begin{align*}
\mathbb{E}\left(\frac{1}{D^{\alpha}}\right)=\sum_{j=1}^{\infty}\frac{p_j}{j^{\alpha}}
=\frac{1}{d}\sum_{j=1}^{\infty}\frac{p_j}{j^{\alpha-1}}=\frac{1}{d}b.
\end{align*}

Then, since
$\mathbb{E}\left(cD^{-1}R^{(k-1)}\right)=c(1-p_0)d^{-1}<\infty$
and $\E(N)=d$, we apply Lemma~\ref{lemma properties}(iii) to
obtain
\begin{align*}
\P&(R^{(k)}>x)\sim
\P\left(\sum_{j=1}^{N}\frac{c}{D_{j}}R^{(k-1)}+[1-c(1-p_0)]>x\right)\\
&\sim\left(c^{\alpha}bC_{k-1}+
\left(\frac{c(1-p_0)}{d}\right)^{\alpha}\right)\P(N>x-[1-c(1-p_0)])\\
&\sim\left(c^{\alpha}bC_{k-1}+
\left(\frac{c(1-p_0)}{d}\right)^{\alpha}\right)\P(N>x)
\;\mbox{as}\;x\to\infty,
\end{align*}
for any $k\ge 2$. Here the last relation again follows from the
property of long-tailed random variables (\ref{eq_long}).{\sloppy

}

Then for the constant $C_k$ we have
\begin{align*}
C_k &= c^{\alpha}\;b\; C_{k-1}+
\left(\frac{c(1-p_0)}{d}\right)^{\alpha}\\
&= \left(c^{\alpha}b
\left(\frac{c(1-p_0)}{d}\right)^{\alpha}\sum_{j=0}^{k-2}c^{j\alpha}b^{j}+
\left(\frac{c(1-p_0)}{d}\right)^{\alpha}\right)\\
&=\left(\frac{c(1-p_0)}{d}\right)^{\alpha}\sum_{j=0}^{k-1}c^{j\alpha}b^{j}.
\end{align*}
\end{proof}

\bibliography{myliterature}

\begin{thebibliography}{10}

\bibitem{IndoChina}
http://law.dsi.unimi.it/.
\newblock Accessed in January 2007.

\bibitem{Stanford}
http://www.stanford.edu/$\sim$sdkamvar/research.html.
\newblock Accessed in March 2006.

\bibitem{Albert99}
R.~Albert and A.~L. Barab\`{a}si.
\newblock Emergence of scaling in random networks.
\newblock {\em Science}, 286:509--512, 1999.

\bibitem{Aldous05}
D.~J. Aldous and A.~Bandyopadhyay.
\newblock A survey of max-type recursive distributional equations.
\newblock {\em Ann. Appl. Probab.}, 15:1047--1110, 2005.

\bibitem{Andersen06}
R.~Andersen, F.~Chung, and K.~Lang.
\newblock Local graph partitioning using pagerank vectors.
\newblock In {\em Proceedings of the 47th Annual IEEE Symposium on Foundations
  of Computer Science (FOCS'06)}, pages 475--486, Washington, DC, USA, 2006.
  IEEE Computer Society.

\bibitem{Avrachenkov06}
K.~Avrachenkov and D.~Lebedev.
\newblock {P}age{R}ank of scale free growing networks.
\newblock Technical Report 5858, INRIA, 2006.

\bibitem{Avrachenkov06a}
K.~Avrachenkov and N.~Litvak.
\newblock The effect of new links on {G}oogle {P}age{R}ank.
\newblock {\em Stoch. Models}, 22(2):319--331, 2006.

\bibitem{Baeza06}
R.~Baeza-Yates, C.~Castillo, and E.~Efthimiadis.
\newblock Characterization of national {W}eb domains.
\newblock {\em To appear in {ACM} {TOIT}}, 2006.

\bibitem{Becchetti06}
L.~Becchetti and C.~Castillo.
\newblock The distribution of {P}age{R}ank follows a power-law only for
  particular values of the damping factor.
\newblock In {\em Proceedings of the 15th international conference on World
  Wide Web}, pages 941--942. ACM Press, New York, 2006.

\bibitem{Becchetti06a}
L.~Becchetti, C.~Castillo, D.~Donato, and A.~Fazzone.
\newblock A comparison of sampling techniques for {W}eb characterization.
\newblock In {\em Workshop on Link Analysis ({LinkKDD})}, 2006.

\bibitem{Bingham_RV}
N.~H. Bingham, C.~M. Goldie, and J.~L. Teugels.
\newblock {\em Regular Variation}.
\newblock Cambridge University Press, 1989.

\bibitem{Brin98}
S.~Brin and L.~Page.
\newblock The anatomy of a large-scale hypertextual web search engine.
\newblock {\em Comput. Networks ISDN Systems}, 33:107--117, 1998.

\bibitem{Broder00}
A.~Broder, R.~Kumar, F.~Maghoul, P.~Raghavan, S.~Rajagopalan, R.~Statac,
  A.~Tomkins, and J.~Wiener.
\newblock Graph structure in the {W}eb.
\newblock {\em Comput. Networks}, 33:309--320, 2000.

\bibitem{Capocci06}
A.~Capocci, V.~D.~P. Servedio, F.~Colaiori, L.~S. Buriol, D.~Donato,
  S.~Leonardiand, and G.~Caldarelli.
\newblock Preferential attachment in the growth of social networks: the case of
  {W}ikipedia.
\newblock Technical Report 0602026, arXiv/physics, 2006.

\bibitem{Chen06}
P.~Chen, H.~Xie, S.~Maslov, and S.~Redner.
\newblock Finding scientific gems with {G}oogle.
\newblock Technical Report 0604130, arxiv/physics/, 2006.

\bibitem{Donato04}
D.~Donato, L.~Laura, S.~Leonardi, and S.~Millozi.
\newblock Large scale properties of the {W}ebgraph.
\newblock {\em Eur. Phys. J.}, 38:239--243, 2004.

\bibitem{Eiron04}
N.~Eiron, K.~S. McCurley, and J.~A. Tomlin.
\newblock Ranking the {W}eb frontier.
\newblock In {\em WWW '04: Proceedings of the 13th international conference on
  World Wide Web}, pages 309--318, New York, NY, USA, 2004. ACM Press.

\bibitem{Fortunato05a}
S.~Fortunato, M.~Boguna, A.~Flammini, and F.Menczer.
\newblock How to make the top ten: Approximating {P}age{R}ank from in-degree,
  2005.
\newblock arXiv.org/cs/cs.IR/0511016.

\bibitem{Fortunato06}
S.~Fortunato and A.~Flammini.
\newblock Random walks on directed networks: the case of {P}age{R}ank, 2006.
\newblock arxiv.org/physics/0604203.

\bibitem{Gyongyi04}
Z.~Gyongyi, H.~Garcia-Molina, and J.~Pedersen.
\newblock Combating web spam with trustrank.
\newblock In {\em 30th International Conference on Very Large Data Bases}, page
  576–587, 2004.

\bibitem{Haveliwala03}
T.H. Haveliwala.
\newblock Topic-sensitive {P}age{R}ank: A context-sensitive ranking algorithm
  for {W}eb search.
\newblock {\em IEEE TKDE}, 15(4):784--796, 2003.

\bibitem{Jessen06}
A.~H. Jessen and T.~Mikosch.
\newblock Regularly varying functions.
\newblock {\em Publications de L'Institut Mathematique, Nouvelle S\'erie},
  79(93), 2006.

\bibitem{Langville03}
A.~N. Langville and C.~D. Meyer.
\newblock Deeper inside {P}age{R}ank.
\newblock {\em Internet Math.}, 1:335--380, 2003.

\bibitem{Litvak06}
N.~Litvak, W.~R.~W. Scheinhardt, and Y.~Volkovich.
\newblock In-degree and {P}age{R}ank: Why do they follow similar power laws?
\newblock {\em \emph{To appear in} Internet Math.}

\bibitem{Liu01}
Q.~Liu.
\newblock Asymptotic properties and absolute continuity of laws stable by
  random weighted mean.
\newblock {\em Stochastic Process. Appl.}, 95(1):83--107, September 2001.

\bibitem{Newman05}
M.~E.~J. Newman.
\newblock Power laws, {P}areto distributions and {Z}ipf's law.
\newblock {\em Contemporary Physics}, 46:323--351, 2005.

\bibitem{Pandurangan02}
G.~Pandurangan, P.~Raghavan, and E.~Upfal.
\newblock {U}sing {P}age{R}ank to characterize web structure.
\newblock In {\em 8th Annual International Computing and Combinatorics
  Conference (COCOON)}, Singapore, 2002.

\bibitem{Ross_SP}
S.~M. Ross.
\newblock {\em Stochastic processes}.
\newblock Wiley Series in Probability and Statistics: Probability and
  Statistics. John Wiley \& Sons Inc., New York, second edition, 1996.

\bibitem{Ross03}
S.~M. Ross.
\newblock The inspection paradox.
\newblock {\em Probab. Engrg. Inform. Sci.}, 17:47--51, 2003.

\end{thebibliography}
\bibliographystyle{plain}
\end{document}